# ZERO BIASING AND A DISCRETE CENTRAL LIMIT THEOREM

By Larry Goldstein and Aihua Xia

*University of Southern California and University of Melbourne*

We introduce a new family of distributions to approximate $\mathbb{P}(W \in A)$ for $A \subset \{\ldots, -2, -1, 0, 1, 2, \ldots\}$ and $W$ a sum of independent integer-valued random variables $\xi_1, \xi_2, \ldots, \xi_n$ with finite second moments, where, with large probability, $W$ is not concentrated on a lattice of span greater than 1. The well-known Berry–Esseen theorem states that, for $Z$ a normal random variable with mean $\mathbb{E}(W)$ and variance $\text{Var}(W)$, $\mathbb{P}(Z \in A)$ provides a good approximation to $\mathbb{P}(W \in A)$ for $A$ of the form $(-\infty, x]$. However, for more general $A$, such as the set of all even numbers, the normal approximation becomes unsatisfactory and it is desirable to have an appropriate discrete, nonnormal distribution which approximates $W$ in total variation, and a discrete version of the Berry–Esseen theorem to bound the error. In this paper, using the concept of zero biasing for discrete random variables (cf. Goldstein and Reinert [*J. Theoret. Probab.* **18** (2005) 237–260]), we introduce a new family of discrete distributions and provide a discrete version of the Berry–Esseen theorem showing how members of the family approximate the distribution of a sum $W$ of integer-valued variables in total variation.

**1. Introduction.** We introduce a new family of distributions to approximate $\mathbb{P}(W \in A)$ for $A$ a subset of $\mathbf{Z} = \{\ldots, -2, -1, 0, 1, 2, \ldots\}$ and $W$ a sum of independent integer-valued random variables $\xi_1, \xi_2, \ldots, \xi_n$ with finite second moments, where the probability that $W$ is not concentrated on a lattice of span greater than 1 is large. When $A$ is of the form $(-\infty, x]$ and $\xi_i$'s have finite third moments, we can use the well-known Berry–Esseen theorem ([7] and [15]) which states that there exists an absolute constant $C$ such that

$$\sup_{z \in \mathbb{R}} \left| \mathbb{P}\left( \frac{W - \mu}{\sigma} \leq z \right) - \Phi(z) \right| \leq \frac{C}{\sigma^3} \sum_{i=1}^{n} \mathbb{E}[|\xi_i - \mathbb{E}(\xi_i)|^3],$$

where $\mu = \mathbb{E}(W)$, $\sigma^2 = \text{Var}(W)$ and $\Phi$ is the cumulative distribution function of the standard normal. If $\xi_i$'s are identically distributed, then the bound is of the order $n^{-1/2}$, which is known to be the best possible. However, for more general $A$, such as the set of all even numbers, the errors of normal approximation may be large, or difficult to compute; for such cases, it is desirable to have a distribution which approximates $W$ in total variation, and a discrete version of the Berry–Esseen theorem to evaluate the error. Moreover, approximations in total variation have the property that any function of $W$ is approximated in total variation to the same degree as $W$ itself, an advantage not enjoyed by the Kolmogorov distance.

A few discrete distributions, such as signed compound Poisson measures and translated Poisson distributions (see [6, 9] and references therein) have been proposed to make very close approximations in total variation to the distribution of $W$. These approximations can be viewed as modifications of Poisson approximation and in applications, one often transforms the sum $W$ into a form which can be approximated reasonably well by a suitably chosen Poisson random variable. In estimating the errors of approximation, besides the assumption that $W$ has large probability of not being concentrated on a lattice of span greater than 1, one also needs other assumptions, such as existence of the third moments of the $\xi_i$'s ([6], Theorem 4.3), and may additionally introduce truncation. Another approach is to define a discrete normal $Y$ by

$$\mathbb{P}(Y = j) = \mathbb{P}(j - 1/2 < Z \leq j + 1/2), \qquad Z \sim \mathcal{N}(\mu, \sigma^2), \qquad j \in \mathbf{Z}$$

(L. H. Y. Chen, personal communication), though it is not clear what quality of approximation $Y$ can achieve.

In this paper we propose a class of approximating distributions which have carrier space $\mathbf{Z}$, thus avoiding truncation and integerization problems. These new distributions are uniquely determined by parameters $\mu$ and $\sigma^2$, similarly to how the approximating normal distribution is determined in the classical central limit theorem. It is expected that any such approximating family of discrete distributions be related to the Poisson, a distribution characterized by the property of being equal to its own reduced Palm distribution; see [23], page 93. As this property is intrinsic in the study of certain Poisson approximations [1, 11], and since the Palm distribution involves only the first moment of the distribution, it is of interest to determine whether there exists any counterpart to the Poisson also involving the second moment, which gives additional flexibility in approximation. One appropriate counterpart can be uncovered through the concept of zero biasing [20]. Based on the continuous normal case, it is expected that the class of approximating distributions should arrive naturally as the unique candidates which equal their zero-biased distribution. However, because of the discrete setting, some adjustments are first needed to make the idea work. In Section 2, we provide some background on zero biasing in both the continuous and discrete settings, and define our approximating family of distributions through a modified zero biasing form. In particular, our distributions are related to the operator (2.11), connected to discrete zero biasing, and are the stationary laws of the processes with corresponding generator (3.3), similarly to how normal laws are related to an operator connected to continuous zero biasing, and are the stationary distributions of the Ornstein–Uhlenbeck processes. Next, in Section 3 we establish the Stein equation and employ the bilateral birth and death processes ([28], Chapter 8) to estimate the Stein factors, in a similar fashion to that in [8]. In Section 4 a general approximation theorem is given which provides a bound in total variation between an integer-valued random variable $Y$ and a member of the family of our approximating distributions, in terms of the distance between $Y$ and its zero-biased distribution, paralleling Lemma 2.1 in [17] for the continuous case. The bound is of the same order as the normal approximation when the weaker Komogorov metric is used. The general theorem is then applied to obtain a bound for the approximation of a sum $W$ of independent integer-valued random variables under (only) second moment conditions, yielding a form which simplifies further under the assumption of finite third moments.

**2. Zero biasing and characterization of the approximating distribution.**
For any nonnegative random variable $X$ with mean $\mathbb{E}(X) = \mu \in (0, \infty)$ and distribution $dF(x)$, the $X$-*size biased* distribution is given by

$$(2.1) \qquad dF^s(x) = \frac{x\,dF(x)}{\mu}, \qquad x \geq 0,$$

or, equivalently, by the characterizing equation

$$\mathbb{E}[Xf(X)] = \mu \mathbb{E} f(X^s) \qquad \text{for all } f \text{ with } \mathbb{E}|Xf(X)| < \infty.$$

It is often helpful to think of size biasing as a transformation defined on nonnegative distributions with finite mean. Size biasing can appear (unwanted and sometimes unnoticed) in various sampling contexts [13]; for example, in random digit dialing, where $F$ in (2.1) is the uniform distribution on telephone numbers, it is twice as likely to dial a household with $x = 2$ telephone lines than a household where $x = 1$. When $X$ is a nonnegative integer-valued random variable with positive finite mean $\mu$, the $X$-size biased distribution (2.1) specializes to

$$(2.2) \qquad \mathbb{P}(X^s = k) = \frac{k\mathbb{P}(X = k)}{\mu}, \qquad k = 0, 1, \ldots.$$

The counterpart of size biasing in point process theory is the Palm distributions (see [23], Chapter 10) introduced by Palm in 1943. It is easily verified that $X$ has a Poisson distribution if and only if $\mathcal{L}(X^s) = \mathcal{L}(X+1)$. This fact can be used to study Poisson approximation and is part of the foundation of the well-known Stein–Chen method (see [10] or [4]).

One notable property of the size-biased transformation is that a sum of independent nonnegative random variables can be size-biased by replacing a single summand, chosen with probability proportional to its mean, with one independent of the remaining variables and having that summand's size-biased distribution; that is, with $\xi_i$ independent nonnegative variables with finite mean $\mathbb{E}\xi_i = \mu_i, i = 1, \ldots, n$, and

$$W = \sum_{i=1}^{n} \xi_i, \qquad \text{we have } W^s = W - \xi_I + \xi_I^s,$$

where $I$ is a random index, independent of $\xi_1, \ldots, \xi_n$, with distribution

$$\mathbb{P}(I = i) = \frac{\mu_i}{\sum_{j=1}^{n} \mu_j}.$$

For $\xi$ a nontrivial indicator variable, (2.2) shows that $\xi^s = 1$. Hence, a sum of independent indicators $\xi_i$, $i = 1, \ldots, n$, can be size-biased by setting a single indicator, chosen with probability proportional to $\mathbb{E}\xi_i$, to one.

The zero bias transformation was introduced in [20], based both on its similarity to the size-biased transformation and the following characterization of the mean zero normal distribution given in [26], which forms the basis of Stein's celebrated method for normal approximation [27]: $Z$ is a mean zero, variance $\sigma^2$ normal variable if and only if, for all absolutely continuous $f$ with $\mathbb{E}|Zf(Z)| < \infty$,

$$\mathbb{E}[Zf(Z)] = \sigma^2 \mathbb{E}f'(Z).$$

For any $Y$ with $\mathbb{E}Y = 0$ and $\text{Var}(Y) = \sigma^2$, Goldstein and Reinert [20] prove that there exists $Y^*$, called the $Y$-*zero biased* distribution, such that, for all absolutely continuous $f$ with $\mathbb{E}|Yf(Y)| < \infty$,

(2.3) $$\mathbb{E}[Yf(Y)] = \sigma^2 \mathbb{E}f'(Y^*).$$

By the Stein (if and only if) characterization, it is clear that $Y$ has the mean zero normal distribution if and only if $\mathcal{L}(Y) = \mathcal{L}(Y^*)$. In other words, the mean zero normal distribution is the unique fixed point of the zero bias transformation. Heuristically, then, one can show that $Y$ is close to normal by showing that $Y$ is close to $Y^*$; for in this case, $Y$ itself is close to being a fixed point and, therefore, should be close to the unique fixed point, the normal. For this reason, it is key that zero biasing enjoys a property similar to the one mentioned above which holds for size biasing. A sum $Y$ of independent mean zero variables $\xi_1, \ldots, \xi_n$ with finite variances $\sigma_1^2, \ldots, \sigma_n^2$ can be zero-biased by choosing a variable using an independent index $I$ with distribution

$$\mathbb{P}(I = i) = \frac{\sigma_i^2}{\sum_{j=1}^n \sigma_j^2},$$

which takes values with probability proportional to variance, and replacing the selected variable with one from that summand's zero-biased distribution which is independent of the remaining variables, that is,

(2.4) $$Y^* = Y - \xi_I + \xi_I^*.$$

Hence, a sum of roughly comparable independent mean zero variables with finite variances is close in distribution to normal, since its zero bias distribution differs from its original one by only one comparable summand of many. For applications of the zero bias transformation to simple random sampling (see [20]) to hierarchical structures (see [17]) to combinatorial central limit theorems [18] and to the computation of $L^1$ bounds to the normal [20].

Goldstein and Reinert [21] show that both size biasing and zero biasing are special cases of distributional transformations specified by a biasing function $P$ and an order $m$; both size biasing and zero biasing have $P(x) = x$, and orders $m = 0$ and $m = 1$, respectively; such transformations are often related to families of orthogonal polynomials. To approximate by a given distribution, one can often construct a transformation for which it is a fixed point. The transformations for which discrete distributions will be fixed points have derivative replaced by difference, in particular, with $\Delta f(i) := f(i+1) - f(i)$, the Poisson distribution with mean $\lambda$ is a fixed point of the transformation characterized by

$$\mathbb{E}[(Y - \lambda)f(Y)] = \lambda \mathbb{E}\Delta f(Y^\star).$$

However, one obtains additional flexibility by not insisting that the mean and variance be equal. Therefore, parallel to (2.3), we give the following definition:

DEFINITION 2.1. For an integer-valued random variable $Y$ with mean $\mu$ and finite variance $\sigma^2$, we say that $Y^\star$ has the discrete $Y$-zero biased distribution if, for all bounded functions $f : \mathbf{Z} \to \mathbf{R}$,

(2.5) $$\mathbb{E}[(Y - \mu)f(Y)] = \sigma^2 \mathbb{E}\Delta f(Y^\star).$$

It is easily verified that (2.4) holds for the discrete zero bias transformation, that is, that a sum of independent discrete random variables can be discrete zero-biased by replacing one variable, chosen with probability proportional to variance, by a variable from that summand's discrete zero bias distribution, independent of the remaining variables. When no confusion between the discrete and continuous cases can arise, we simply say that $Y^\star$ has the $Y$-zero biased distribution.

For $Y$ an integer-valued random variable with finite mean and variance, the existence and uniqueness of $\mathcal{L}(Y^\star)$ can be proved as follows. For each $j \in \mathbf{Z}$, let $f_j(i) = \mathbf{1}_{[j,\infty)}(i)$, so that

$$\Delta f_j(i) = \begin{cases} 1, & \text{for } i = j-1, \\ 0, & \text{for } i \neq j-1, \end{cases}$$

giving

$$\begin{aligned}
\mathbb{P}(Y^\star = j - 1) &= \frac{\mathbb{E}[(Y - \mu)\mathbf{1}_{(Y \geq j)}]}{\sigma^2} \\
&= \frac{\mathbb{E}[Y \mathbf{1}_{[j,\infty)}(Y)] - \mu \mathbb{P}(Y \geq j)}{\sigma^2}.
\end{aligned} \quad (2.6)$$

For $j \geq \mu$, (2.6) is clearly nonnegative, and the identity $\mathbb{E}[(Y - \mu)\mathbf{1}_{(Y \geq j)}] = -\mathbb{E}[(Y - \mu)\mathbf{1}_{(Y < j)}]$ implies (2.6) is also nonnegative for $j < \mu$. Using this identity and the fact that $\sum_{j=1}^\infty \mathbf{1}_{(Y \geq j)} = Y\mathbf{1}_{(Y \geq 1)}$ and $\sum_{j=-\infty}^0 \mathbf{1}_{(Y < j)} = -Y\mathbf{1}_{(Y \leq 0)}$, we have the probabilities in (2.6) summing to one, since

$$\begin{aligned}
\mathbb{E}\left[\sum_{j=1}^\infty (Y - \mu)\mathbf{1}_{(Y \geq j)} + \sum_{j=-\infty}^0 (Y - \mu)\mathbf{1}_{(Y \geq j)}\right] \\
= \mathbb{E}\left[\sum_{j=1}^\infty (Y - \mu)\mathbf{1}_{(Y \geq j)} - \sum_{j=-\infty}^0 (Y - \mu)\mathbf{1}_{(Y < j)}\right] \\
= \mathbb{E}[(Y - \mu)Y\mathbf{1}_{(Y \geq 1)} + (Y - \mu)Y\mathbf{1}_{(Y \leq 0)}] \\
= \mathbb{E}[(Y - \mu)Y] = \sigma^2.
\end{aligned}$$

For $\eta$ an indicator variable with $\text{Var}(\eta) = \theta^2 = (1 - \mathbb{E}\eta)\mathbb{E}\eta > 0$, (2.6) shows that $\eta^\star = 0$:

$$\begin{aligned}
\mathbb{P}(\eta^\star = 0) &= \frac{\mathbb{E}[(\eta - \mathbb{E}\eta)\mathbf{1}_{(\eta \geq 1)}]}{\theta^2} = \frac{\mathbb{E}[(\eta - \mathbb{E}\eta)\mathbf{1}_{(\eta = 1)}]}{\theta^2} \\
&= \frac{(1 - \mathbb{E}\eta)\mathbb{P}(\eta = 1)}{\theta^2} = 1.
\end{aligned} \quad (2.7)$$

Though true in this particular case, it is incorrect to conclude from this example that $\eta^\star = \eta^s - 1$, that is, that the discrete zero bias operation is the same as the reduced Palm. For an independent sum, the Palm distribution is obtained by replacing a summand chosen proportional to its mean, but to achieve the zero bias distribution, one chooses proportional to variance.

Since the fixed points of the continuous zero bias transformation (2.3) are the mean zero normal distributions, it is of immediate interest to determine which distributions, if any, are fixed points of the discrete zero bias transformation (2.5), that is, to find which $\tilde{S}$ satisfy

$$(2.8) \qquad \mathbb{E}[(\tilde{S} - \mu)f(\tilde{S})] = \sigma^2 \mathbb{E}\Delta f(\tilde{S})$$

for all bounded functions $f$ on $\mathbf{Z}$. We show now that, unlike the situation in the continuous case, distributional fixed points do *not* exist for all choices of $\mu$ and $\sigma^2$. It is for this reason that we introduce the family of distributions given in Lemma 2.2. Using $f_j(i) = \mathbf{1}_{(i=j)}$ for $j \in \mathbf{Z}$ in (2.8), we have

$$(2.9) \qquad (\sigma^2 + j - \mu)\mathbb{P}(\tilde{S} = j) = \sigma^2 \mathbb{P}(\tilde{S} = j - 1), \qquad j \in \mathbf{Z}.$$

There are two cases to check for (2.9), depending on whether or not $\mu - \sigma^2$ is an integer. Let $\tilde{\kappa} := \min\{i : i \geq \mu - \sigma^2\}$. When $\mu - \sigma^2$ is an integer, (2.9) gives that $\mathbb{P}(\tilde{S} = j) = 0$ for $j < \tilde{\kappa}$. However, when $\mu - \sigma^2$ is not an integer, then, unless $\mathcal{L}(\tilde{S})$ is the null measure on $\mathbf{Z}$, the values $\mathbb{P}(\tilde{S} = \tilde{\kappa} - 1)$ and $\mathbb{P}(\tilde{S} = \tilde{\kappa} - 2)$ strictly alternate in sign, that is, $\mathcal{L}(\tilde{S})$ is a signed measure which takes on both positive and negative values. To avoid such a signed measure when $\mu - \sigma^2$ is not an integer, we truncate the distribution at, for example, $\tilde{\kappa}$, so that $\mathbb{P}(\tilde{S} = j) = 0$ for $j < \tilde{\kappa}$, in which case (2.9) fails for $j = \tilde{\kappa}$ and $\tilde{S}$ is only approximately a fixed point of the discrete zero bias transformation. In either of these cases, where $\mu - \sigma^2$ is an integer or where $\mu - \sigma^2$ is not an integer and we truncate at $\tilde{\kappa}$, iteration of (2.9) yields

$$\mathbb{P}(\tilde{S} = \tilde{\kappa}) = \left\{ \sum_{j=\tilde{\kappa}+1}^{\infty} \left( \prod_{\tilde{\kappa}+1}^{j} \frac{\sigma^2}{\sigma^2 + i - \mu} \right) + 1 \right\}^{-1},$$

$$\mathbb{P}(\tilde{S} = j) = \left( \prod_{i=\tilde{\kappa}+1}^{j} \frac{\sigma^2}{\sigma^2 + i - \mu} \right) \mathbb{P}(\tilde{S} = \tilde{\kappa}), \qquad j \geq \tilde{\kappa} + 1.$$

If $\mu - \sigma^2$ is an integer we now see that $\tilde{S}$ corresponds to a translated Poisson ([6], page 131); that is, the distribution of $\tilde{S}$ equals that of $Y + \mu - \sigma^2$ with $Y$ a Poisson random variable with mean $\sigma^2$. Further, elementary calculations using (2.9) yield

$$(2.10) \qquad \mathbb{E}\tilde{S} = \mu \mathbb{P}(\tilde{S} \geq \tilde{\kappa} + 1) + (\sigma^2 + \tilde{\kappa})\mathbb{P}(\tilde{S} = \tilde{\kappa}),$$

so $\mathbb{E}\tilde{S} = \mu$ if and only if $\mu - \sigma^2$ is an integer; since (2.10) will not be used later on, we omit the details.

Using a truncated approximating distribution, such as $\tilde{S}$ above, leaves $\mathbb{P}(W < \tilde{\kappa})$ in the upper bound when we estimate the error caused by approximating $W$ by $\tilde{S}$, and can become quite inconvenient in applications ([4], Section 9.2 and [2]). To avoid truncation, we introduce a two-parameter family of distributions which have carrier space $\mathbf{Z}$, similar to the two-parameter normal distributions which have carrier space $\mathbb{R}$. For an integer $\kappa$ with $\mu - \sigma^2 \leq \kappa < \mu + \sigma^2 + 1$, define the operator

$$(2.11) \quad \mathcal{B}f(i) = \begin{cases} \sigma^2 f(i+1) - (\sigma^2 + i - \mu)f(i) \\ \quad = \sigma^2 \Delta f(i) - (i - \mu)f(i), & i \geq \kappa, \\ (\sigma^2 + \mu - i)f(i+1) - \sigma^2 f(i) \\ \quad = \sigma^2 \Delta f(i) - (i - \mu)f(i+1), & i \leq \kappa - 1, \end{cases}$$

for all bounded functions $f$ on $\mathbf{Z}$. Note that $\sigma^2 + i - \mu$ and $\sigma^2 + \mu - i$ are nonnegative over their respective ranges $i \geq \kappa$ and $i \leq \kappa - 1$, and strictly positive except when $\mu - \sigma^2$ is an integer and $i = \kappa = \mu - \sigma^2$. The following lemmas are devoted to the properties of $\mathcal{L}(S) := \Psi_\kappa(\mu, \sigma^2)$, the distribution characterized by $\mathcal{B}$:

LEMMA 2.2. *There exists a unique distribution $\mathcal{L}(S) = \Psi_\kappa(\mu, \sigma^2)$, characterized by $\mathbb{E}\mathcal{B}f(S) = 0$ for all bounded functions $f$ on $\mathbf{Z}$, whose distribution $\pi_i = \mathbb{P}(S = i)$, $i \in \mathbf{Z}$, satisfies*

$$(2.12) \quad \pi_\kappa = \left\{ \sum_{j=\kappa+1}^{\infty} \left( \prod_{\kappa+1}^{j} \frac{\sigma^2}{\sigma^2 + i - \mu} \right) + 1 + \frac{\sigma^2 + \kappa - \mu}{\sigma^2 + \mu - \kappa + 1} + \sum_{j=-\infty}^{\kappa-2} \left( \prod_{i=j}^{\kappa-2} \frac{\sigma^2}{\sigma^2 + \mu - i} \right) \frac{\sigma^2 + \kappa - \mu}{\sigma^2 + \mu - \kappa + 1} \right\}^{-1}$$

*and*

$$\pi_j = \begin{cases} \left( \prod_{i=\kappa+1}^{j} \frac{\sigma^2}{\sigma^2 + i - \mu} \right) \pi_\kappa, & j \geq \kappa + 1, \\ \dfrac{\sigma^2 + \kappa - \mu}{\sigma^2 + \mu - \kappa + 1} \pi_\kappa, & j = \kappa - 1, \\ \left( \prod_{i=j}^{\kappa-2} \frac{\sigma^2}{\sigma^2 + \mu - i} \right) \left( \dfrac{\sigma^2 + \kappa - \mu}{\sigma^2 + \mu - \kappa + 1} \right) \pi_\kappa, & j \leq \kappa - 2. \end{cases}$$

*Moreover, $\mathbb{E}(|S|^l) < \infty$ for $0 \leq l < \infty$ and $\mathbb{E}\mathcal{B}f(S) = 0$ for all $f$ such that $\mathbb{E}\{|S|[|f(S)| + |f(S+1)|]\} < \infty$.*

REMARK. When $\mu - \sigma^2$ is an integer and $\kappa = \mu - \sigma^2$, the distribution of $S$ reduces to that of $\tilde{S}$ when $\tilde{\kappa} = \kappa$.

PROOF OF LEMMA 2.2. Since, for each fixed $j$ and $\mathbf{1}_j(i) = \mathbf{1}_{(i=j)}$,

$$0 = \mathbb{E}[\mathcal{B}\mathbf{1}_j(S)] = \sum_i \mathcal{B}\mathbf{1}_j(i)\pi_i = \mathcal{B}\mathbf{1}_j(j-1)\pi_{j-1} + \mathcal{B}\mathbf{1}_j(j)\pi_j,$$

we obtain recursive formulæ as follows:

$$-(\sigma^2 + j - \mu)\pi_j + \sigma^2 \pi_{j-1} = 0, \qquad j \geq \kappa + 1, \tag{2.13}$$

$$-(\sigma^2 + j - \mu)\pi_j + (\sigma^2 + \mu - j + 1)\pi_{j-1} = 0, \qquad j = \kappa, \tag{2.14}$$

$$-\sigma^2 \pi_j + (\sigma^2 + \mu - j + 1)\pi_{j-1} = 0, \qquad j \leq \kappa - 1. \tag{2.15}$$

Hence,

$$\pi_j = \left( \prod_{i=\kappa+1}^{j} \frac{\sigma^2}{\sigma^2 + i - \mu} \right) \pi_\kappa, \qquad j \geq \kappa + 1,$$

$$\pi_{\kappa-1} = \frac{\sigma^2 + \kappa - \mu}{\sigma^2 + \mu - \kappa + 1} \pi_\kappa \tag{2.16}$$

and

$$\pi_{j-1} = \frac{\sigma^2}{\sigma^2 + \mu - j + 1} \pi_j = \left( \prod_{i=j-1}^{\kappa-2} \frac{\sigma^2}{\sigma^2 + \mu - i} \right) \pi_{\kappa-1}, \qquad j \leq \kappa - 1,$$

so, replacing $j$ by $j+1$ in the last identity, it follows from (2.16) that

$$\pi_j = \left( \prod_{i=j}^{\kappa-2} \frac{\sigma^2}{\sigma^2 + \mu - i} \right) \left( \frac{\sigma^2 + \kappa - \mu}{\sigma^2 + \mu - \kappa + 1} \right) \pi_\kappa, \qquad j \leq \kappa - 2.$$

Summing the probabilities to one yields (2.12). Convergence is guaranteed, for the sum in (2.12) over $j \geq \kappa + 1$, say, by the fact that $\sigma^2/(\sigma^2 + i - \mu) \leq \sigma^2/(i - \kappa)$ for all $i \geq \kappa + 1$ and the fact that $\sum_{j=\kappa+1}^{\infty} \sigma^{2(j-\kappa)}/(j - \kappa)! < \infty$. Hence, the distribution of $S$ exists and is uniquely determined by the specified distribution.

The claim $\mathbb{E}(|S|^l) < \infty$ follows from the fact that

$$\sum_{j \geq \kappa+1} |j|^l \pi_j \leq \pi_\kappa \sum_{j \geq \kappa+1} |j|^l \frac{\sigma^{2(j-\kappa)}}{(j - \kappa)!} < \infty$$

and

$$\sum_{j \leq \kappa-2} |j|^l \pi_j \leq \pi_{\kappa-1} \sum_{j \leq \kappa-2} |j|^l \frac{\sigma^{2(\kappa-j-1)}}{(\kappa - j - 1)!} < \infty.$$

Finally, taking $f_n = (f \wedge n) \vee (-n)$, $n = 1, 2, \ldots$, we have $\mathbb{E}\mathcal{B}f_n(S) = 0$ and

$$|\mathcal{B}f_n(i)| \leq (|i| + |\mu| + 2\sigma^2)[|f(i)| + |f(i+1)|].$$

Hence, the dominated convergence theorem ensures that $\mathbb{E}\mathcal{B}f(S) = 0$ by letting $n \to \infty$. □

LEMMA 2.3.  $\mathbb{E}(S) = \mu$ and $\mathrm{Var}(S) = \sigma^2 + (\sigma^2 + \kappa - \mu)\pi_\kappa$.

PROOF. Letting $f(i) \equiv 1$, since $\mathcal{B}f(i) = \mu - i$ for all $i$, $\mathbb{E}\mathcal{B}f(S) = 0$ yields $\mathbb{E}S = \mu$. Next, letting $f(i) = i$ in (2.11), for $i \geq \kappa$ we have $\mathcal{B}f(i) = \sigma^2 - i^2 + \mu i$, while for $i \leq \kappa - 1$ we have $\mathcal{B}f(i) = \sigma^2 - i^2 + \mu i + \mu - i$, which can be written as

$$\mathcal{B}f(i) = \sigma^2 - i^2 + \mu i + (\mu - i)\mathbf{1}_{(i \leq \kappa - 1)}.$$

It follows from Lemma 2.2 that $\mathbb{E}\mathcal{B}f(S) = 0$, which yields

$$0 = \sigma^2 - \text{Var}(S) + \mathbb{E}(\mu - S)\mathbf{1}_{(S \leq \kappa - 1)}$$

so that $\text{Var}(S) = \sigma^2 + \mathbb{E}(\mu - S)\mathbf{1}_{(S \leq \kappa - 1)}$.

Now, using $\mathbb{E}(S - \mu) = 0$, and (2.13) for the fourth equality, it follows that

$$\mathbb{E}(\mu - S)\mathbf{1}_{(S \leq \kappa - 1)} = \mathbb{E}(S - \mu)\mathbf{1}_{(S \geq \kappa)} = \sum_{i \geq \kappa}(i - \mu)\pi_i$$

$$= \sum_{i \geq \kappa + 1}(\sigma^2 + i - \mu)\pi_i - \sigma^2 \sum_{i \geq \kappa + 1}\pi_i + (\kappa - \mu)\pi_\kappa$$

$$= \sum_{i \geq \kappa + 1}\sigma^2 \pi_{i-1} - \sigma^2 \sum_{i \geq \kappa + 1}\pi_i + (\kappa - \mu)\pi_\kappa$$

$$= (\sigma^2 + \kappa - \mu)\pi_\kappa.$$

Hence,

$$\text{Var}(S) = \sigma^2 + (\sigma^2 + \kappa - \mu)\pi_\kappa. \qquad \square$$

Note that if we choose $\kappa = \min\{i : i \geq \mu - \sigma^2\}$, then $|\text{Var}(S) - \sigma^2| \leq \pi_\kappa$. The following lemma shows in what sense $S$ is close to a fixed point of the zero bias transformation when $\text{Var}(S)$ is close to $\sigma^2$:

LEMMA 2.4. *The S-zero biased distribution $S^\star$, given in Definition 2.1, satisfies*

$$\mathbb{P}(S^\star = j) = \begin{cases} \sigma^2 \mathbb{P}(S = j)/\text{Var}(S), & j \geq \kappa, \\ 1 - \sigma^2/\text{Var}(S), & j = \kappa - 1, \\ \sigma^2 \mathbb{P}(S = j + 1)/\text{Var}(S), & j \leq \kappa - 2. \end{cases}$$

PROOF. Fixing $j \geq \kappa$ and letting $f_j(i) = \mathbf{1}_{[j+1,\infty)}(i)$, we have $\mathcal{B}f_j(i) = \sigma^2 \Delta f_j(i) - (i - \mu)f_j(i), \forall i \in \mathbf{Z}$. Using the characterization equation $\mathbb{E}\mathcal{B}f_j(S) = 0$ and Definition 2.1,

$$0 = \mathbb{E}(\sigma^2 \Delta f_j(S) - (S - \mu)f_j(S)) = \mathbb{E}(\sigma^2 \Delta f_j(S) - \text{Var}(S)\Delta f_j(S^\star))$$

which, along with $\Delta f_j(i) = \mathbf{1}_{(i=j)}$, gives the claim for $j \geq \kappa$.

Likewise, fixing $j \leq \kappa - 2$ and letting $f_j(i) = \mathbf{1}_{(-\infty, j+1]}(i)$, we have $\mathcal{B}f_j(i) = \sigma^2 \Delta f_j(i) - (i - \mu)f_j(i+1), \forall i \in \mathbf{Z}$, and

$$0 = \mathbb{E}(\sigma^2 \Delta f_j(S) - (S - \mu)f_j(S+1))$$
$$= \mathbb{E}(\sigma^2 \Delta f_j(S) - \mathrm{Var}(S)\Delta f_j(S^\star + 1)),$$

which, with $\Delta f_j(i) = -\mathbf{1}_{(i=j+1)}$, gives the claim for $j \leq \kappa - 2$. Finally, the value $\mathbb{P}(S^\star = \kappa - 1)$ can be obtained from $\sum_{i=-\infty}^{\infty} \mathbb{P}(S^\star = i) = 1$. □

**3. Stein's method and Stein's factors.** Brown and Xia [8] introduced a class of approximating distributions $\pi$, determined by parameters $\alpha_i, \beta_i, i \in \mathbf{Z}_+ := \{0, 1, 2, \ldots\}$, satisfying

(3.1) $$\pi_i \alpha_i = \pi_{i+1} \beta_{i+1}, \qquad i \in \mathbf{Z}_+.$$

Equation (3.1) enabled the authors of that work to view $\pi$ as the stationary distribution of a birth–death process and to give a neat probabilistic derivation of *Stein magic factors*, essentially under the condition that for each $k = 1, 2, \ldots$,

(3.2) $$\alpha_k - \alpha_{k-1} \leq \beta_k - \beta_{k-1},$$

letting $\beta_0 = 0$. A key point in that derivation is that the solution to the Stein equation is an explicit linear combination of mean upward and downward transition times of the birth–death process ([8], Lemma 2.1). Under condition (3.2), all differences of the solution of the Stein equation are negative except one—an essential structure for the neat derivation of Stein magic factors for polynomial birth–death approximations, which include Poisson, binomial, negative binomial and hypergeometric approximations [8].

In this section we consider approximating distributions $\pi$ on $\mathbf{Z}$ (instead of $\mathbf{Z}_+$) which are determined by two parameters $\mu$ and $\sigma^2$ and which satisfy the balance equation (3.4). Analogously to the context in [8], we define a generator (3.3) of a bilateral birth–death process such that $\pi$ is its stationary distribution. In Lemma 3.4, we prove that all differences of the solution of the Stein equation are negative except one and derive the Stein magic factors.

For each bounded function $g$ on $\mathbf{Z}$, writing $f(x+1) = g(x+1) - g(x)$, we have

(3.3)
$$\mathcal{B}f(i) = \begin{cases} \sigma^2(g(i+1) - g(i)) \\ \quad + (\sigma^2 + i - \mu)(g(i-1) - g(i)), & i \geq \kappa, \\ (\sigma^2 + \mu - i)(g(i+1) - g(i)) \\ \quad + \sigma^2(g(i-1) - g(i)), & i \leq \kappa - 1, \end{cases}$$
$$:= \mathcal{A}g(i) = \alpha_i(g(i+1) - g(i)) + \beta_i(g(i-1) - g(i)),$$

where

$$\alpha_i = \begin{cases} \sigma^2, & i \geq \kappa, \\ \sigma^2 + \mu - i, & i \leq \kappa - 1, \end{cases}$$

and

$$\beta_i = \begin{cases} \sigma^2 + i - \mu, & i \geq \kappa, \\ \sigma^2, & i \leq \kappa - 1. \end{cases}$$

$\mathcal{A}$ is the generator of bilateral birth and death processes ([28], Chapter 8) with "birth rates" specified by $\{\alpha_i : i \in \mathbf{Z}\}$ and "death rates" by $\{\beta_i : i \in \mathbf{Z}\}$. When $\mu - \sigma^2 < \kappa < \mu + \sigma^2 + 1$ so that all $\alpha_i$'s and $\beta_i$'s are positive, the bilateral birth and death processes are always nonexplosive and ergodic ([28], Chapter 8 and [12]). However, when $\mu - \sigma^2$ is an integer and $\kappa = \mu - \sigma^2$, we get $\beta_\kappa = 0$ (the only possible zero of all the transition rates), which means that all states in $(-\infty, \kappa - 1]$ are transient. In other words, when the Markov chain is at a state in $(-\infty, \kappa - 1]$, it will move quickly into states in $[\kappa, \infty)$, while if the Markov chain is at a state in $[\kappa, \infty)$, it will never visit states in $(-\infty, \kappa - 1]$. In this case, the approximating distribution is the same as a translated Poisson; it has been well treated in various papers (see [6, 9] and references therein).

Čekanavičius and Vaitkus [9] studied the translated Poisson (referred to as *centered Poisson* in the paper) approximation to the sum $W$ of independent indicator random variables with $\lambda = \mathbb{E}(W)$ and $\lambda_2 = \lambda - \text{Var}(W)$. Their approximating translated Poisson is the sum of $\lfloor \lambda_2 \rfloor$, the integer part of $\lambda_2$, and a Poisson random variable with mean $\lambda - \lfloor \lambda_2 \rfloor$. This distribution is a slight variation of our $\tilde{S}$, and a straightforward modification of the Stein–Chen method is used to estimate the approximation errors. Hence, from now on, we concentrate on the case where $\mu - \sigma^2 < \kappa < \mu + \sigma^2 + 1$.

It is a routine exercise to check that $\Psi_\kappa(\mu, \sigma^2)$ is the equilibrium distribution of the Markov chain with generator $\mathcal{A}$, and that it satisfies the following balance equation:

$$\alpha_i \pi_i = \beta_{i+1} \pi_{i+1} \qquad \forall i \in \mathbf{Z}. \tag{3.4}$$

Denote by $Z_i(t), t \geq 0$, the Markov chain generated by $\mathcal{A}$ with initial value $i$, and define stopping times

$$\begin{aligned}
\tau_i &= \inf\{t : Z_i(t) \neq i\}, \\
\tau_i^+ &= \inf\{t : Z_i(t) = i + 1\}, \\
\tau_i^- &= \inf\{t : Z_i(t) = i - 1\}, \qquad i \in \mathbf{Z}.
\end{aligned} \tag{3.5}$$

LEMMA 3.1. *For every bounded function $h$ on $\mathbf{Z}$, the integral*

$$g(i) := -\int_0^\infty \{\mathbb{E}[h(Z_i(t))] - \mathbb{E}[h(S)]\} \, dt$$

*is well defined and satisfies the* Stein *identity*

$$\mathcal{A}g(i) = h(i) - \mathbb{E}h(S).$$

PROOF. Split the bilateral birth–death process $Z_i$ at $\kappa$ into two ordinary birth–death processes. Each of the two processes is a standard linear model, hence exponentially ergodic, implying that the process $Z_i$ is also exponentially ergodic (see [12], Theorem 4.1, or [14], page 1679). More precisely, taking $\tilde{V}(i) = 1 + |i - \kappa|$, $c_0 = 1$ and $b_0$ sufficiently large, we can see that the condition $(\tilde{D})$ (see the remark after the statement of the condition) in page 1679 of [14] is satisfied, meaning that, by the second paragraph of [14], page 1681, there is some $0 < \rho < 1$ such that, for all $i \in \mathbf{Z}$, there exists a finite constant $M_i$ with

$$\sum_{j \in \mathbf{Z}} |\mathbb{P}(Z_i(t) = j) - \pi_j| \leq M_i \rho^t \qquad \text{for all } t \geq 0.$$

Hence,

$$\int_0^\infty |\mathbb{E}[h(Z_i(t))] - \mathbb{E}[h(S)]| \, dt \leq \sup_{j \in \mathbf{Z}} |h(j)| \int_0^\infty M_i \rho^t \, dt < \infty,$$

which ensures that $g$ is well defined.

Next, the general theory of Markov processes ensures that, for $a > 0$,

$$(a - \mathcal{A})^{-1}(h - \mathbb{E}h(S))(i) = \int_0^\infty e^{-at}[\mathbb{E}h(Z_i(t)) - \mathbb{E}h(S)] \, dt,$$

(see [16], page 10), and the Stein identity corresponds to the above equation when $a = 0$. A sketch of the proof of the Stein identity is as follows. Since $\tau_i \sim \exp(\alpha_i + \beta_i)$ and

$$\mathbb{P}(\tau_i = \tau_i^+) = \frac{\alpha_i}{\alpha_i + \beta_i} \quad \text{and} \quad \mathbb{P}(\tau_i = \tau_i^-) = \frac{\beta_i}{\alpha_i + \beta_i},$$

by invoking the strong Markov property and momentarily ignoring integrability issues, we get

$$\begin{aligned}
g(i) &= -\mathbb{E} \int_0^\infty [h(Z_i(t)) - \mathbb{E}h(S)] \, dt \\
&= -\mathbb{E}\left\{ \int_0^{\tau_i} [h(Z_i(t)) - \mathbb{E}h(S)] \, dt + \int_{\tau_i}^\infty [h(Z_i(t)) - \mathbb{E}h(S)] \, dt \right\} \\
&= -\frac{h(i) - \mathbb{E}h(S)}{\alpha_i + \beta_i} - \mathbb{E} \int_0^\infty [h(Z_i(t + \tau_i)) - \mathbb{E}h(S)] \, dt \\
&= -\frac{h(i) - \mathbb{E}h(S)}{\alpha_i + \beta_i} - \frac{\alpha_i}{\alpha_i + \beta_i} \mathbb{E} \int_0^\infty [h(Z_{i+1}(t)) - \mathbb{E}h(S)] \, dt \\
&\quad - \frac{\beta_i}{\alpha_i + \beta_i} \mathbb{E} \int_0^\infty [h(Z_{i-1}(t)) - \mathbb{E}h(S)] \, dt \\
&= -\frac{h(i) - \mathbb{E}h(S)}{\alpha_i + \beta_i} + \frac{\alpha_i}{\alpha_i + \beta_i} g(i+1) + \frac{\beta_i}{\alpha_i + \beta_i} g(i-1),
\end{aligned}$$

(3.6)

which, after reorganizing the terms, implies

$$h(i) - \mathbb{E}h(S) = \alpha_i(g(i+1) - g(i)) + \beta_i(g(i-1) - g(i)) = \mathcal{A}g(i),$$

as desired. To prove (3.6) rigorously, by the strong Markov property, we have, for each $0 < u < \infty$,

$$-\int_0^u \mathbb{E}[h(Z_i(t)) - h(S)]\,dt$$

$$= -\mathbb{E}\int_0^u [h(Z_i(t)) - \mathbb{E}h(S)]\,dt$$

$$= -\mathbb{E}\int_0^{\tau_i \wedge u} [h(Z_i(t)) - \mathbb{E}h(S)]\,dt - \mathbb{E}\int_{\tau_i \wedge u}^u [h(Z_i(t)) - \mathbb{E}h(S)]\,dt$$

$$= [\mathbb{E}h(S) - h(i)]\mathbb{E}(\tau_i \wedge u)$$
$$\quad - \int_0^u \mathbb{E}\left\{\int_s^u [h(Z_i(t)) - \mathbb{E}h(S)]\,dt \,\Big|\, \tau_i = s\right\}\mathbb{P}(\tau_i \in ds)$$

$$= [\mathbb{E}h(S) - h(i)]\mathbb{E}(\tau_i \wedge u)$$
$$\quad - \int_0^u \mathbb{E}\left\{\int_0^{u-s} [h(Z_i(s+v)) - \mathbb{E}h(S)]\,dv \,\Big|\, \tau_i = s\right\}\mathbb{P}(\tau_i \in ds)$$

$$= [\mathbb{E}h(S) - h(i)]\mathbb{E}(\tau_i \wedge u)$$
$$\quad - \frac{\alpha_i}{\alpha_i + \beta_i}\int_0^u \left\{\int_0^{u-s} \mathbb{E}[h(Z_{i+1}(v)) - h(S)]\,dv\right\}\mathbb{P}(\tau_i \in ds)$$
$$\quad - \frac{\beta_i}{\alpha_i + \beta_i}\int_0^u \left\{\int_0^{u-s} \mathbb{E}[h(Z_{i-1}(v)) - h(S)]\,dv\right\}\mathbb{P}(\tau_i \in ds)$$

$$= [\mathbb{E}h(S) - h(i)]\mathbb{E}(\tau_i \wedge u)$$
$$\quad - \frac{\alpha_i}{\alpha_i + \beta_i}\int_0^u \mathbb{E}[h(Z_{i+1}(v)) - h(S)]\mathbb{P}(\tau_i \leq u - v)\,dv$$
$$\quad - \frac{\beta_i}{\alpha_i + \beta_i}\int_0^u \mathbb{E}[h(Z_{i-1}(v)) - h(S)]\mathbb{P}(\tau_i \leq u - v)\,dv.$$

Letting $u \to \infty$ and applying the bounded convergence theorem yields (3.6).
□

For fixed $k_1,\ k_2 \in \mathbf{Z}$ with $k_1 \leq k_2$, define

$$e_i^-(k_1, k_2) = \mathbb{E}\int_0^{\tau_i^-} \mathbf{1}_{[k_1, k_2]}(Z_i(t))\,dt$$

and

$$e_i^+(k_1, k_2) = \mathbb{E}\int_0^{\tau_i^+} \mathbf{1}_{[k_1, k_2]}(Z_i(t))\,dt,$$

the expected time that the Markov chain $Z_i(t)$ spends in $[k_1, k_2]$ before it reaches $i-1$ and $i+1$, respectively. We note that $e_i^-(-\infty, \infty) = \mathbb{E}\tau_i^- := \overline{\tau_i^-}$ and $e_i^+(-\infty, \infty) = \mathbb{E}\tau_i^+ := \overline{\tau_i^+}$, as introduced in [8], page 1378. Hence, the following lemma generalizes Lemma 2.2 in [8]:

LEMMA 3.2.
$$e_i^-(k_1, k_2) = \begin{cases} \frac{\sum_{l=i\vee k_1}^{k_2} \pi_l}{\beta_i \pi_i}, & \text{if } i \leq k_2, \\ 0, & \text{if } i > k_2, \end{cases}$$

and

$$e_i^+(k_1, k_2) = \begin{cases} \frac{\sum_{l=k_1}^{i\wedge k_2} \pi_l}{\alpha_i \pi_i}, & \text{if } i \geq k_1, \\ 0, & \text{if } i < k_1. \end{cases}$$

PROOF. Since $\tau_i \sim \exp(\alpha_i + \beta_i)$, and $\tau_i \leq \tau_i^-$ by (3.5), we have

$$e_i^-(k_1, k_2) = \mathbb{E}\int_0^{\tau_i} \mathbf{1}_{[k_1,k_2]}(Z_i(t))\,dt + \mathbb{E}\int_{\tau_i}^{\tau_i^-} \mathbf{1}_{[k_1,k_2]}(Z_i(t))\,dt$$

$$= \mathbf{1}_{\{k_1 \leq i \leq k_2\}} \frac{1}{\alpha_i + \beta_i}$$

$$+ \mathbb{E}\left(\int_{\tau_i}^{\tau_i^-} \mathbf{1}_{[k_1,k_2]}(Z_i(t))\,dt \Big| \tau_i = \tau_i^-\right)\mathbb{P}(\tau_i = \tau_i^-)$$

$$+ \mathbb{E}\left(\int_{\tau_i}^{\tau_i^-} \mathbf{1}_{[k_1,k_2]}(Z_i(t))\,dt \Big| \tau_i < \tau_i^-\right)\mathbb{P}(\tau_i < \tau_i^-).$$

The second-to-last term is clearly zero. For the last term, given $\tau_i < \tau_i^-$, we have $Z_i(\tau_i) = i+1$, so by the strong Markov property,

$$\mathbb{E}\left(\int_{\tau_i}^{\tau_i^-} \mathbf{1}_{[k_1,k_2]}(Z_i(t))\,dt \Big| \tau_i < \tau_i^-\right) = \mathbb{E}\int_0^{\tau_{i+1,i-1}} \mathbf{1}_{[k_1,k_2]}(Z_{i+1}(t))\,dt,$$

where $\tau_{j_1, j_2} = \inf\{t \colon Z_{j_1}(t) = j_2\}$. Now, again, by the strong Markov property,

$$\mathbb{E}\int_0^{\tau_{i+1,i-1}} \mathbf{1}_{[k_1,k_2]}(Z_{i+1}(t))\,dt$$

$$= \mathbb{E}\int_0^{\tau_{i+1}^-} \mathbf{1}_{[k_1,k_2]}(Z_{i+1}(t))\,dt + \mathbb{E}\int_{\tau_{i+1}^-}^{\tau_{i+1,i-1}} \mathbf{1}_{[k_1,k_2]}(Z_{i+1}(t))\,dt$$

$$= \mathbb{E}\int_0^{\tau_{i+1}^-} \mathbf{1}_{[k_1,k_2]}(Z_{i+1}(t))\,dt + \mathbb{E}\int_0^{\tau_i^-} \mathbf{1}_{[k_1,k_2]}(Z_i(t))\,dt.$$

Combining the equations above gives

$$e_i^-(k_1,k_2) = \mathbf{1}_{\{k_1 \leq i \leq k_2\}} \frac{1}{\alpha_i + \beta_i} + (e_{i+1}^-(k_1,k_2) + e_i^-(k_1,k_2))\frac{\alpha_i}{\alpha_i + \beta_i},$$

which, using (3.4), implies

$$\pi_i \beta_i e_i^-(k_1,k_2) = \pi_i \mathbf{1}_{\{k_1 \leq i \leq k_2\}} + \pi_i \alpha_i e_{i+1}^-(k_1,k_2)$$
$$= \pi_i \mathbf{1}_{\{k_1 \leq i \leq k_2\}} + \beta_{i+1} \pi_{i+1} e_{i+1}^-(k_1,k_2).$$

Clearly, $e_i^-(k_1,k_2) = 0$ for $i > k_2$, so

$$\pi_i \beta_i e_i^-(k_1,k_2) = \sum_{l=i}^{k_2} [\pi_l \beta_l e_l^-(k_1,k_2) - \pi_{l+1} \beta_{l+1} e_{l+1}^-(k_1,k_2)] = \sum_{l=i}^{k_2} \pi_l \mathbf{1}_{\{k_1 \leq l \leq k_2\}},$$

which implies

$$e_i^-(k_1,k_2) = \begin{cases} \dfrac{\sum_{l=i}^{k_2} \pi_l}{\beta_i \pi_i}, & \text{if } k_1 \leq i \leq k_2, \\ \dfrac{\sum_{l=k_1}^{k_2} \pi_l}{\beta_i \pi_i}, & \text{if } i < k_1, \end{cases}$$

as desired.

Likewise,

$$e_i^+(k_1,k_2) = \mathbb{E}\int_0^{\tau_i} \mathbf{1}_{[k_1,k_2]}(Z_i(t))\,dt + \mathbb{E}\int_{\tau_i}^{\tau_i^+} \mathbf{1}_{[k_1,k_2]}(Z_i(t))\,dt$$
$$= \mathbf{1}_{\{k_1 \leq i \leq k_2\}} \frac{1}{\alpha_i + \beta_i} + (e_{i-1}^+(k_1,k_2) + e_i^+(k_1,k_2))\frac{\beta_i}{\alpha_i + \beta_i},$$

which, together with (3.4), gives

$$\pi_i \alpha_i e_i^+(k_1,k_2) = \pi_i \mathbf{1}_{\{k_1 \leq i \leq k_2\}} + \pi_i \beta_i e_{i-1}^+(k_1,k_2)$$
$$= \pi_i \mathbf{1}_{\{k_1 \leq i \leq k_2\}} + \alpha_{i-1} \pi_{i-1} e_{i-1}^+(k_1,k_2).$$

We have that $e_i^+(k_1,k_2) = 0$ for $i < k_1$, so

$$\pi_i \alpha_i e_i^+(k_1,k_2) = \sum_{l=k_1}^{i} [\pi_l \alpha_l e_l^+(k_1,k_2) - \alpha_{l-1} \pi_{l-1} e_{l-1}^+(k_1,k_2)] = \sum_{l=k_1}^{i} \pi_l \mathbf{1}_{\{k_1 \leq l \leq k_2\}},$$

again giving the claimed expression. □

Note that in the sequel, we will only need the quantities $\mathbb{E}\tau_i^- = e_i^-(-\infty,\infty)$ and $\mathbb{E}\tau_i^+ = e_i^+(-\infty,\infty)$, since we will focus on the choice $\kappa = \min\{i : i \geq \mu\}$ and the total variation metric. For other cases, the general result in Lemma 3.2 is needed.

Since
$$\mathbb{P}(W \in A) - \mathbb{P}(S \in A) = \sum_{j \in A}[\mathbb{P}(W = j) - \mathbb{P}(S = j)]$$
$$= \sum_{j \in A}[\mathbb{E}\mathbf{1}_{\{j\}}(W) - \pi_j],$$

we define
$$h_j(x) = \mathbf{1}_{\{j\}}(x) - \pi_j,$$

$g_j$ to be the solution of $\mathcal{A}g_j = h_j$, and $f_j(i) = g_j(i) - g_j(i-1)$.

LEMMA 3.3. *For each $j \in \mathbf{Z}$,*
$$f_j(i) = \begin{cases} -\pi_j \dfrac{\sum_{l=-\infty}^{i-1} \pi_l}{\alpha_{i-1}\pi_{i-1}}, & \text{for } i \leq j, \\ \pi_j \dfrac{\sum_{l=i}^{\infty} \pi_l}{\beta_i \pi_i}, & \text{for } i > j. \end{cases}$$

PROOF.  Using the strong Markov property, we have
$$\int_0^u \mathbb{E}h_j(Z_{i-1}(t))\,dt = \mathbb{E}\int_0^{u \wedge \tau_{i-1}^+} h_j(Z_{i-1}(t))\,dt + \mathbb{E}\int_{u \wedge \tau_{i-1}^+}^u h_j(Z_{i-1}(t))\,dt$$
$$= \mathbb{E}\int_0^{u \wedge \tau_{i-1}^+} h_j(Z_{i-1}(t))\,dt$$
$$+ \int_0^u \left\{\int_0^{u-s} \mathbb{E}h_j(Z_i(v))\,dv\right\}\mathbb{P}(\tau_{i-1}^+ \in ds),$$

and letting $u \to \infty$ yields
$$\int_0^\infty \mathbb{E}h_j(Z_{i-1}(t))\,dt = \mathbb{E}\int_0^{\tau_{i-1}^+} h_j(Z_{i-1}(t))\,dt + \int_0^\infty \mathbb{E}h_j(Z_i(t))\,dt.$$

Hence, for $i \leq j$,
$$f_j(i) = g_j(i) - g_j(i-1)$$
$$= -\int_0^\infty \mathbb{E}h_j(Z_i(t))\,dt + \int_0^\infty \mathbb{E}h_j(Z_{i-1}(t))\,dt$$
$$= \mathbb{E}\int_0^{\tau_{i-1}^+} h_j(Z_{i-1}(t))\,dt$$
$$= -\pi_j \mathbb{E}\tau_{i-1}^+ = -\pi_j \frac{\sum_{l=-\infty}^{i-1} \pi_l}{\alpha_{i-1}\pi_{i-1}},$$

the last equality following from Lemma 3.2. Likewise, using

$$\int_0^\infty \mathbb{E} h_j(Z_i(t))\,dt = \mathbb{E}\int_0^{\tau_i^-} h_j(Z_i(t))\,dt + \int_0^\infty \mathbb{E} h_j(Z_{i-1}(t))\,dt$$

and Lemma 3.2 again, it follows that, for $i > j$,

$$f_j(i) = -\int_0^\infty \mathbb{E} h_j(Z_i(t))\,dt + \int_0^\infty \mathbb{E} h_j(Z_{i-1}(t))\,dt$$

$$= -\mathbb{E}\int_0^{\tau_i^-} h_j(Z_i(t))\,dt = \pi_j \mathbb{E}\tau_i^- = \pi_j \frac{\sum_{l=i}^\infty \pi_l}{\beta_i \pi_i}. \qquad \square$$

LEMMA 3.4. *For $A \subset \mathbf{Z}$, let*

$$h_A(x) = \mathbf{1}_A(x) - \mathbb{P}(S \in A),$$

*$g_A$ be the solution to*

$$\mathcal{A} g_A = h_A \quad \text{and} \quad f_A(i) = g_A(i) - g_A(i-1).$$

*If $\kappa = \min\{i : i \geq \mu\}$, then for all $i$ and $A$,*

$$|\Delta f_A(i)| \leq \frac{1-\pi_i}{\alpha_i \wedge \beta_i} \wedge \frac{1}{\alpha_i} \wedge \frac{1}{\beta_i} \leq \frac{1-\pi_i}{\sigma^2}.$$

REMARK. For approximating distributions on $\mathbf{Z}_+$ satisfying the balance equation (3.1) ([8], page 1382) proved that, if (3.2) is satisfied, then $|\Delta f_A(i)| \leq \frac{1}{\alpha_i} \wedge \frac{1}{\beta_i}$, and [22], Corollary 3.5.1, gives the bound $|\Delta f_A(i)| \leq \Delta f_i(i)$ under the assumption of nonincreasing $\alpha_i$'s and nondecreasing $\beta_i$'s, derived similarly to the inequality $\Delta f_i(i) \geq \Delta f_j(i)$ below. Lemma 3.4 is parallel to these types of estimates for the new approximating distribution satisfying the version of condition (3.2) which has been appropriately modified for its range.

PROOF OF LEMMA 3.4. It follows from Lemma 3.3 that

$$\Delta f_j(i) = f_j(i+1) - f_j(i)$$

$$= \begin{cases} -\pi_j \left( \dfrac{\sum_{l=-\infty}^i \pi_l}{\alpha_i \pi_i} - \dfrac{\sum_{l=-\infty}^{i-1} \pi_l}{\alpha_{i-1} \pi_{i-1}} \right), & i < j, \\ \pi_j \left( \dfrac{\sum_{l=i+1}^\infty \pi_l}{\beta_{i+1} \pi_{i+1}} + \dfrac{\sum_{l=-\infty}^{i-1} \pi_l}{\alpha_{i-1} \pi_{i-1}} \right), & i = j, \\ \pi_j \left( \dfrac{\sum_{l=i+1}^\infty \pi_l}{\beta_{i+1} \pi_{i+1}} - \dfrac{\sum_{l=i}^\infty \pi_l}{\beta_i \pi_i} \right), & i > j. \end{cases}$$

Since $\kappa = \min\{i : i \geq \mu\}$ and, therefore, $\mu \leq \kappa \leq \mu + 1$, one can verify directly that $\{\alpha_i,\ i \in \mathbf{Z}\}$ are nonincreasing and $\{\beta_i,\ i \in \mathbf{Z}\}$ are nondecreasing. Hence, for $i < j$,

$$\frac{\sum_{l=-\infty}^{i} \pi_l}{\alpha_i \pi_i} - \frac{\sum_{l=-\infty}^{i-1} \pi_l}{\alpha_{i-1} \pi_{i-1}} = \frac{\sum_{l=-\infty}^{i} \pi_l}{\alpha_i \pi_i} - \frac{\sum_{l=-\infty}^{i-1} \pi_l}{\beta_i \pi_i}$$

$$= \frac{1}{\alpha_i \beta_i \pi_i} \left( \beta_i \sum_{l=-\infty}^{i} \pi_l - \alpha_i \sum_{l=-\infty}^{i-1} \pi_l \right)$$

$$\geq \frac{1}{\alpha_i \beta_i \pi_i} \left( \sum_{l=-\infty}^{i} \beta_l \pi_l - \sum_{l=-\infty}^{i-1} \alpha_l \pi_l \right) = 0,$$

where for the first and last equalities we have applied the balance equation $\alpha_l \pi_l = \beta_{l+1} \pi_{l+1}$. Likewise, for $i > j$,

$$\frac{\sum_{l=i+1}^{\infty} \pi_l}{\beta_{i+1} \pi_{i+1}} - \frac{\sum_{l=i}^{\infty} \pi_l}{\beta_i \pi_i} = \frac{\sum_{l=i+1}^{\infty} \pi_l}{\alpha_i \pi_i} - \frac{\sum_{l=i}^{\infty} \pi_l}{\beta_i \pi_i}$$

$$= \frac{1}{\alpha_i \beta_i \pi_i} \left( \beta_i \sum_{l=i+1}^{\infty} \pi_l - \alpha_i \sum_{l=i}^{\infty} \pi_l \right)$$

$$\leq \frac{1}{\alpha_i \beta_i \pi_i} \left( \sum_{l=i+1}^{\infty} \beta_l \pi_l - \sum_{l=i}^{\infty} \alpha_l \pi_l \right) = 0.$$

Hence, $\Delta f_j(i) \leq 0$ for $j \neq i$ and $\Delta f_j(i) > 0$ for $j = i$, and for any $A \subset \mathbf{Z}$,

$$\Delta f_A(i) = \sum_{j \in A} \Delta f_j(i) \leq \Delta f_i(i) = \pi_i \left( \frac{\sum_{l=i+1}^{\infty} \pi_l}{\alpha_i \pi_i} + \frac{\sum_{l=-\infty}^{i-1} \pi_l}{\beta_i \pi_i} \right)$$

(3.7)
$$= \frac{\sum_{l=i+1}^{\infty} \pi_l}{\alpha_i} + \frac{\sum_{l=-\infty}^{i-1} \pi_l}{\beta_i}$$

$$\leq \frac{1}{\alpha_i \wedge \beta_i} \left( \sum_{l=i+1}^{\infty} \pi_l + \sum_{l=-\infty}^{i-1} \pi_l \right) = \frac{1 - \pi_i}{\alpha_i \wedge \beta_i}.$$

To obtain the other terms in the bound, note that, since $\{\alpha_i, i \in \mathbf{Z}\}$ are nonincreasing and $\{\beta_i, i \in \mathbf{Z}\}$ are nondecreasing, for $l \geq i+1$, we have

$$\pi_l = \frac{\alpha_{l-1} \pi_{l-1}}{\beta_l} \leq \frac{\alpha_i \pi_{l-1}}{\beta_i},$$

and for $l \leq i - 1$,

$$\pi_l = \frac{\beta_{l+1} \pi_{l+1}}{\alpha_l} \leq \frac{\beta_i \pi_{l+1}}{\alpha_i},$$

which in turn imply

(3.8)
$$\frac{\sum_{l=i+1}^{\infty} \pi_l}{\alpha_i} \leq \frac{\sum_{l=i+1}^{\infty} \pi_{l-1}}{\beta_i} = \frac{1}{\beta_i} \sum_{l=i}^{\infty} \pi_l$$

and

$$\frac{\sum_{l=-\infty}^{i-1} \pi_l}{\beta_i} \leq \frac{\sum_{l=-\infty}^{i-1} \pi_{l+1}}{\alpha_i} = \frac{\sum_{l=-\infty}^{i} \pi_l}{\alpha_i}. \quad (3.9)$$

Now, it follows from (3.7) and (3.8) that $\Delta f_A(i) \leq 1/\beta_i$, while combining (3.7) and (3.9) gives $\Delta f_A(i) \leq 1/\alpha_i$.

On the other hand, since $h_\mathbf{Z} = 1 - \mathbb{P}(S \in \mathbf{Z}) \equiv 0$, we have

$$-\Delta f_A(i) = \Delta f_{\mathbf{Z}\setminus A}(i) \leq \frac{1-\pi_i}{\alpha_i \wedge \beta_i} \wedge \frac{1}{\alpha_i} \wedge \frac{1}{\beta_i}.$$

Noting that $\alpha_i$ and $\beta_i$ are both at least $\sigma^2$ for all $i$ completes the proof. □

**4. Zero biasing and approximation theorems.** We define the total variation distance between two probability measures $\mathbf{Q}_1, \mathbf{Q}_2$ on $\mathbf{Z}$ as follows:

$$d_{\mathrm{TV}}(\mathbf{Q}_1, \mathbf{Q}_2) = \sup_{A \subset \mathbf{Z}} |\mathbf{Q}_1(A) - \mathbf{Q}_2(A)|.$$

Using Lemma 3.4, we prove the following general theorems, which parallel results in the continuous case showing that $Y$ is close to normal when $Y$ and $Y^*$ are close. Throughout this section, we write $\Psi(\mu, \sigma^2)$ for $\Psi_\kappa(\mu, \sigma^2)$ for $\kappa$ chosen as in Section 3, that is,

$$\kappa = \min\{i : i \geq \mu\}.$$

THEOREM 4.1. *Let $Y$ be an integer-valued random variable with mean $\mu$ and finite variance $\sigma^2$, and let $Y^\star$ have the $Y$-zero biased distribution. Then*

$$d_{\mathrm{TV}}(\mathcal{L}(Y), \Psi(\mu, \sigma^2)) \leq \sum_{i=\kappa}^{\infty} |\mathbb{P}(Y = i) - \mathbb{P}(Y^\star = i)|$$
$$+ \sum_{i=-\infty}^{\kappa-1} |\mathbb{P}(Y = i) - \mathbb{P}(Y^\star + 1 = i)|.$$

PROOF. With $h_A$ and $f_A$ as in Lemma 3.4, recalling the form of the operator $\mathcal{B}$ in (2.11), we have, by the zero bias property (2.5),

$$|\mathbb{P}(Y \in A) - \mathbb{P}(S \in A)|$$
$$= |\mathbb{E}h_A(Y)| = |\mathbb{E}\mathcal{A}g_A(Y)| = |\mathbb{E}\mathcal{B}f_A(Y)|$$
$$= |\sigma^2 \mathbb{E}\Delta f_A(Y) - \mathbb{E}\{(Y-\mu)f_A(Y)\mathbf{1}_{Y \geq \kappa} + (Y-\mu)f_A(Y+1)\mathbf{1}_{Y \leq \kappa-1}\}|$$
$$= \sigma^2 |\mathbb{E}\{\Delta f_A(Y)\} - \mathbb{E}\{[\Delta(f_A(Y^\star)\mathbf{1}_{Y^\star \geq \kappa})] + [\Delta(f_A(Y^\star+1)\mathbf{1}_{Y^\star \leq \kappa-1})]\}|.$$

However, note that, for any $\rho$,

$$\Delta(f(i)\mathbf{1}_{i\geq\rho}) = f(i+1)\mathbf{1}_{i+1\geq\rho} - f(i)\mathbf{1}_{i\geq\rho}$$
$$= [\Delta f(i)]\mathbf{1}_{i\geq\rho} + f(\rho)\mathbf{1}_{i=\rho-1}$$

and

$$\Delta(f(i+1)\mathbf{1}_{i\leq\rho-1}) = f(i+2)\mathbf{1}_{i+1\leq\rho-1} - f(i+1)\mathbf{1}_{i\leq\rho-1}$$
$$= [\Delta f(i+1)]\mathbf{1}_{i\leq\rho-2} - f(\rho)\mathbf{1}_{i=\rho-1}.$$

Hence, with the help of the cancellation of the term $f(\rho)\mathbf{1}_{i=\rho-1}$, the above expectation equals

$$\sigma^2|\mathbb{E}\{\Delta f_A(Y)\} - \mathbb{E}\{[\Delta f_A(Y^\star)]\mathbf{1}_{Y^\star \geq \kappa} + [\Delta f_A(Y^\star+1)]\mathbf{1}_{Y^\star \leq \kappa-2}\}|$$

$$= \sigma^2 \left| \sum_{i=-\infty}^{\infty} \Delta f_A(i)\mathbb{P}(Y=i) \right.$$

$$\left. - \sum_{i=\kappa}^{\infty} \Delta f_A(i)\mathbb{P}(Y^\star=i) - \sum_{i=-\infty}^{\kappa-1} \Delta f_A(i)\mathbb{P}(Y^\star+1=i) \right|$$

$$\leq \sigma^2 \sum_{i=\kappa}^{\infty} |\Delta f_A(i)(\mathbb{P}(Y=i) - \mathbb{P}(Y^\star=i))|$$

$$+ \sigma^2 \sum_{i=-\infty}^{\kappa-1} |\Delta f_A(i)(\mathbb{P}(Y=i) - \mathbb{P}(Y^\star+1=i))|$$

$$\leq \sum_{i=\kappa}^{\infty} |\mathbb{P}(Y=i) - \mathbb{P}(Y^\star=i)| + \sum_{i=-\infty}^{\kappa-1} |\mathbb{P}(Y=i) - \mathbb{P}(Y^\star+1=i)|,$$

where we have applied the bound $|\Delta f_A(i)| \leq 1/\sigma^2$ shown in Lemma 3.4. □

Before applying Theorem 4.1 to the case where $W$ is a sum, we note that the existence of a finite first moment of $Y^\star$ is equivalent to the existence of a finite third moment of $Y$; letting $f(y) = y^2$ in (2.5),

$$\mathbb{E}[Y^3 - \mu Y^2] = \sigma^2 \mathbb{E}[2Y^\star + 1].$$

THEOREM 4.2. *Let $\xi_i$, $i = 1, \ldots, n$, be independent integer-valued random variables and let $W = \sum_{i=1}^{n} \xi_i$. Then, with $W_i = W - \xi_i$, $\mathrm{Var}(\xi_i) = \sigma_i^2$ and $\xi_i^\star$ defined on the same space as $\xi_i$, with the $\xi_i$ zero-biased distribution for $i = 1, \ldots, n$, with $\mu = \mathbb{E}(W)$ and $\sigma^2 = \mathrm{Var}(W)$, we have for any $K > 0$,*

$$d_{\mathrm{TV}}(\mathcal{L}(W), \Psi(\mu, \sigma^2))$$

$$\leq \frac{2}{\sigma^2} \sum_{i=1}^{n} \sigma_i^2 d_+^{(i)}[\mathbb{E}(|\xi_i - \xi_i^\star| \wedge K) + \mathbb{E}(|\xi_i - (\xi_i^\star + 1)| \wedge K)]$$

$$+ \frac{2}{\sigma^2} \sum_{i=1}^{n} \sigma_i^2 \left\{ \sum_{|k_1 - k_2| > K} \mathbb{P}(\xi_i = k_1, \xi_i^\star = k_2) \right.$$

$$+ \sum_{|k_1-k_2|>K} \mathbb{P}(\xi_i = k_1, \xi_i^\star + 1 = k_2)\bigg\},$$

where $d_+^{(i)} = d_{\mathrm{TV}}(\mathcal{L}(W_i), \mathcal{L}(W_i+1))$, $i = 1, \ldots, n$. In particular, letting $K \uparrow \infty$,

$$(4.1) \quad d_{\mathrm{TV}}(\mathcal{L}(W), \Psi(\mu, \sigma^2)) \leq \frac{2}{\sigma^2} \sum_{i=1}^n \sigma_i^2 d_+^{(i)} [\mathbb{E}|\xi_i - \xi_i^\star| + \mathbb{E}|\xi_i - (\xi_i^\star + 1)|],$$

which will be finite when $\mathbb{E}|\xi_i|^3 < \infty, i = 1, \ldots, n$.

PROOF. Considering the first sum in the bound of Theorem 4.1, by invoking (2.4) we have

$$\sum_{j=\kappa}^\infty |\mathbb{P}(W=j) - \mathbb{P}(W^\star = j)|$$

$$\leq \sum_{i=1}^n \frac{\sigma_i^2}{\sigma^2} \sum_{j=\kappa}^\infty |\mathbb{P}(W_i + \xi_i = j) - \mathbb{P}(W_i + \xi_i^\star = j)|$$

$$\leq \sum_{i=1}^n \frac{\sigma_i^2}{\sigma^2} \sum_{j=\kappa}^\infty \sum_{k_1,k_2} |\mathbb{P}(W_i = j+k_1) - \mathbb{P}(W_i = j+k_2)| \mathbb{P}(\xi_i = -k_1, \xi_i^\star = -k_2)$$

$$\leq \sum_{i=1}^n \frac{\sigma_i^2}{\sigma^2} \sum_{j=\kappa}^\infty \sum_{|k_1-k_2|\leq K} \sum_{l=k_1 \wedge k_2}^{k_1 \vee k_2 - 1} |\mathbb{P}(W_i = j+l)$$

$$- \mathbb{P}(W_i = j + l + 1)| \mathbb{P}(\xi_i = -k_1, \xi_i^\star = -k_2)$$

$$+ \sum_{i=1}^n \frac{\sigma_i^2}{\sigma^2} \sum_{j=\kappa}^\infty \sum_{|k_1-k_2|>K} [\mathbb{P}(W_i = j+k_1)$$

$$+ \mathbb{P}(W_i = j+k_2)] \mathbb{P}(\xi_i = -k_1, \xi_i^\star = -k_2)$$

$$\leq 2 \sum_{i=1}^n \frac{\sigma_i^2}{\sigma^2} \bigg\{ d_+^{(i)} \sum_{|k_1-k_2|\leq K} |k_1 - k_2| \mathbb{P}(\xi_i = -k_1, \xi_i^\star = -k_2)$$

$$+ \sum_{|k_1-k_2|>K} \mathbb{P}(\xi_i = -k_1, \xi_i^\star = -k_2)\bigg\}$$

$$\leq \frac{2}{\sigma^2} \sum_{i=1}^n \sigma_i^2 \bigg\{ d_+^{(i)} \mathbb{E}(|\xi_i - \xi_i^\star| \wedge K) + \sum_{|k_1-k_2|>K} \mathbb{P}(\xi_i = -k_1, \xi_i^\star = -k_2)\bigg\};$$

the bound on the remaining sum can be shown similarly. □

REMARK. When $W$ is the sum of many terms of comparable order, the bound in (4.1) is small when $d_+^{(i)}$, $i = 1, \ldots, n$, are small, which is ensured by the condition that, with large probability, $W$ is not concentrated on a lattice of span greater than 1; see Remark 4.5.

REMARK. Note that no signed measures, truncation or translation are required, in contrast to Barbour and Xia [6], Čekanavičius and Vaitkus [9] and Barbour and Choi [3].

COROLLARY 4.3. *Let $I_i$, $i = 1, \ldots, n$, be independent indicator random variables with*

$$\mathbb{P}(I_i = 1) = 1 - \mathbb{P}(I_i = 0) = p_i, \qquad i = 1, \ldots, n,$$

$$W = \sum_{i=1}^n I_i, \qquad \mu = \sum_{i=1}^n p_i, \qquad \sigma^2 = \sum_{i=1}^n p_i(1 - p_i) \quad and$$

$$\vartheta^2 = \sigma^2 - \max_i p_i(1 - p_i).$$

*Then*

$$d_{\mathrm{TV}}(\mathcal{L}(W), \Psi(\mu, \sigma^2)) \leq \frac{1}{\vartheta}.$$

As for approximations using the central limit theorem, we do not expect the $p_i$'s to be small; the bound here has the same order as those in the classical central limit theorem, polynomial birth–death approximation [8] and compound Poisson signed measures approximation [6]. Moreover, there are no additional assumptions required as in [8] or signed measures as in [6].

PROOF OF COROLLARY 4.3. Since $W_i$ is unimodal in this case ([25], page 1273), we have

$$d_{\mathrm{TV}}(\mathcal{L}(W_i), \mathcal{L}(W_i + 1)) \leq \max_j \mathbb{P}(W_i = j) \leq \frac{1}{2}(\sigma^2 - p_i(1 - p_i))^{-1/2} \leq \frac{1}{2\vartheta},$$

where the second inequality is due to Barbour and Jensen [5], page 78. Since $I_i^\star = 0$ by (2.7), we have

$$\mathbb{E}|I_i - I_i^\star| = p_i, \qquad \mathbb{E}|I_i - (I_i^\star + 1)| = 1 - p_i,$$

and it follows from (4.1) that

$$d_{\mathrm{TV}}(\mathcal{L}(W), \Psi(\mu, \sigma^2)) \leq \frac{1}{\vartheta \sigma^2} \sum_{i=1}^n \sigma_i^2 [\mathbb{E}|I_i - I_i^\star| + \mathbb{E}|I_i - (I_i^\star + 1)|]$$

$$= \frac{1}{\vartheta}. \qquad \square$$

REMARK. Note that the proofs do not depend on the order of the index set $\{1,\ldots,n\}$ of the $\xi_i$'s, so one may apply the approximation theorems to the sum of independent integer-valued random variables on an arbitrary index set.

To estimate $d_+^{(i)}$ in general, one may apply Proposition 4.6 of [6], quoted below.

PROPOSITION 4.4. *Suppose that $\xi_i$, $1 \leq i \leq n$, are independent integer-valued random variables, and set $u_i = 1 - d_{\mathrm{TV}}(\mathcal{L}(\xi_i), \mathcal{L}(\xi_i + 1))$, $U = \sum_{i=1}^n \min\{u_i, 1/2\}$. Then, if $W = \sum_{i=1}^n \xi_i$, we have*

$$d_{\mathrm{TV}}(\mathcal{L}(W), \mathcal{L}(W+1)) \leq U^{-1/2}.$$

*Hence, with $W_i = W - \xi_i$,*

$$\max_{1 \leq i \leq n} d_{\mathrm{TV}}(\mathcal{L}(W_i), \mathcal{L}(W_i + 1)) \leq (U-1)^{-1/2}.$$

REMARK. As discussed in [24], Section II.12–14, $d_{\mathrm{TV}}(\mathcal{L}(W), \mathcal{L}(W+1))$ is of order $n^{-1/2}$ when $\xi_i$, $i = 1, \ldots, n$, are independent and identically distributed with an aperiodic distribution.

REMARK 4.5. The assumption of aperiodicity is essential here, where the total variation metric is used. To see why, take $\xi_i, i = 1, \ldots, n$, independent with distribution $\mathbb{P}(\xi_i = 0) = \mathbb{P}(\xi_i = 3) = 1/2$. Then, with probability one, $W$ is concentrated on $\{0, 3, 6, \ldots\}$, a lattice of span greater than 1, and

$$\begin{aligned}&d_{\mathrm{TV}}(\mathcal{L}(W), \Psi(\mu, \sigma^2)) \\ &= \tfrac{1}{2} \sum_j |\mathbb{P}(W=j) - \mathbb{P}(S=j)| \geq \tfrac{1}{2} \sum_{j \notin \{0,3,\ldots\}} \mathbb{P}(S=j) = O(1).\end{aligned}$$

If one wants to lift the assumption of aperiodicity, it is essential to weaken the metric to the Kolmogorov metric, in which case, unless higher moments of $\xi_i$'s (e.g., the third moments) do not exist, the Berry–Esseen theorem would be sufficient.

REMARK. When $\kappa = \min\{i : i \geq \mu\}$, the variance of $S$ does not match that of the sum $W$ of $n$ independent and identically distributed integer-valued random variables; however, crude estimates show that $\mathrm{Var}(S)/\mathrm{Var}(W)$ approaches 1 as $n \to \infty$. It is hoped that future research could address this issue and sharpen the estimates of the approximation errors.

**Acknowledgments.** We would like to thank Professor Mufa Chen for providing us the information about the exponential ergodicity of the bilateral birth–death processes. Part of this work was written while A. Xia was visiting the Institute of Mathematics, Academia Sinica, Taiwan, and A. Xia would like to thank Chii-Ruey Hwang for splendid hospitality. L. Goldstein would like to thank the Department of Mathematics and Statistics at the University of Melbourne for a most enjoyable and productive visit.

DEPARTMENT OF MATHEMATICS KAP 108
UNIVERSITY OF SOUTHERN CALIFORNIA
LOS ANGELES, CALIFORNIA 90089-2532
USA
E-MAIL: larry@paradox.usc.edu

DEPARTMENT OF MATHEMATICS AND STATISTICS
UNIVERSITY OF MELBOURNE
VICTORIA 3010
AUSTRALIA
E-MAIL: xia@ms.unimelb.edu.au